\documentclass[a4paper,11pt,reqno,twosided]{amsart} 
\usepackage{amsmath,amssymb,amsthm}
\usepackage{graphpap,latexsym,epsf}
\usepackage{epsf}
\usepackage{color,graphics}
\usepackage[numbers,sort&compress,longnamesfirst]{natbib}

\theoremstyle{plain}
\newtheorem{theorem}{Theorem}[section]  
\newtheorem{corollary}[theorem]{Corollary}  
\newtheorem{lemma}[theorem]{Lemma}  

\newtheorem*{theorem*}{Theorem}

\theoremstyle{plain}  
\newtheorem{remark}[theorem]{Remark}

\newtheoremstyle{citing}
  {3pt}
  {3pt}
  {\itshape}
  {}
  {\bfseries}
  {.}
  {.5em}
  {\thmnote{#3}}

\theoremstyle{citing}

\numberwithin{equation}{section}

\newcommand{\rand}{\partial} 
\newcommand{\where}{\,|\,}

\newcommand{\laplace}{\Delta}

\newcommand{\di}{\;d}  
  
\newcommand{\nz}{{\mathbb N}}

\newcommand{\rz}{{\mathbb R}}

\newcommand{\eps}{\varepsilon}  
\renewcommand{\phi}{\varphi}

\renewcommand{\div }{{\rm div}\,}
\renewcommand{\a }{\alpha }

\newcommand{\Di }{\mathcal{D}^{1,2}_{a,\l}(\R^N) } 
\newcommand{\e }{\varepsilon }

\renewcommand{\l }{\lambda }

\newcommand{\n }{\nabla }

\newcommand{\R}{\mathbb{R}} 
 
\newcommand{\N}{\mathbb{N}}

\begin{document}
 
\title[]
{Perturbation results of critical elliptic equations\\ of
Caffarelli-Kohn-Nirenberg type
}

\author{Veronica Felli}
\author{Matthias Schneider}
\thanks{ V. F. is supported by M.U.R.S.T. under the national project ``Variational Methods and Nonlinear Differential Equations'' and M. S. research is supported by a S.I.S.S.A. postdoctoral fellowship.}
\address{Scuola Internazionale di Studi Avanzati\\
S.I.S.S.A.\\
Via Beirut 2-4\\
34014 Trieste, Italy}
\email{schneid@sissa.it, felli@sissa.it}

\date{\today}  
\begin{abstract}
We find for small $\eps$ positive solutions to the equation 
\[-\div (|x|^{-2a}\n u)-\displaystyle{\frac{\l}{|x|^{2(1+a)}}}\,u= 
\Big(1+\eps k(x)\Big)\frac{u^{p-1}}{|x|^{bp}}\]
in $\rz^N$, which branch off from
the manifold of minimizers in the class of radial functions of the
corresponding Caffarelli-Kohn-Nirenberg type inequality. Moreover, 
our analysis highlights the symmetry-breaking phenomenon in these
inequalities, namely the
existence of non-radial minimizers.  
\end{abstract}

\maketitle

\section{Introduction}
We will consider the following elliptic equation in $\R^N$ in dimension
$N\geq3$
\begin{align}\label{eq:eq1}
-\div (|x|^{-2a}\n u)-\displaystyle{\frac{\l}{|x|^{2(1+a)}}}\,u=
 K(x)\frac{u^{p-1}}{|x|^{bp}},\quad x\in \R^N \setminus\{0\}
\end{align}
where
\begin{align}
\label{eq:6}
\begin{split}
&-\infty<a<\frac{N-2}2,\quad  -\infty<
\l<\left(\frac{N-2a-2}2\right)^2\\
&p=p(a,b)=\displaystyle{\frac{2N}{N-2(1+a-b)}}\quad\mbox{and}\quad a\leq b< a+1.  
\end{split}
\end{align}
For $\lambda=0$ equation (\ref{eq:eq1}) is related to a family of inequalities given by \citet{CKN}, 
\begin{align}\label{eq:CKN}
\|u\|_{p,b}^2:= \left(\int_{\R^N}|x|^{-bp}|u|^p\,dx\right)^{2/p}\leq {\mathcal
  C}_{a,b}\int_{\R^N}|x|^{-2a}|\n u|^2\,dx \qquad \forall u \in C_0^\infty(\rz^N).
\end{align}
For sharp constants and extremal functions we refer to \citet{CatrinaWang}.\\
The natural functional space to study (\ref{eq:eq1}) is 
$D_a^{1,2}(\rz^N)$ defined as the completion of $C^{\infty}_0(\R^N)$ with respect
to the norm
\[
\|\nabla u\|_a:= \|u\|_*=\left[\int_{\R^N}|x|^{-2a}|\nabla
  u|^2\,dx\right]^{1/2}.
\]
We will mainly deal with the perturbative case $K(x)=1+\e k(x)$, namely with
the problem
\begin{equation*}\tag{${\mathcal P}_{a,b,\l}$}
\begin{cases}
-\div (|x|^{-2a}\n u)-\displaystyle{\frac{\l}{|x|^{2(1+a)}}}\,u=
 \big(1+\e k(x)\big)\frac{u^{p-1}}{|x|^{bp}}\\
u\in D_a^{1,2}(\rz^N)   ,\quad u>0 \text{ in } \rz^N\backslash\{0\}.
\end{cases}
\end{equation*}
Concerning the perturbation $k$ we assume
\begin{align}\label{eq:k}
k\in L^{\infty}(\R^N)\cap C(\rz^N).
\end{align}
Our approach is based on an abstract perturbative variational method discussed by
\citet{AmBa1}, which splits our procedure in three main steps. First we consider the unperturbed
problem, i.e. $\eps=0$, and find a one dimensional manifold of radial solutions.
If this manifold is non-degenerate (see Theorem \ref{t:nondeg} below) a one dimensional
reduction of the perturbed variational problem in $D_a^{1,2}(\rz^N)$ is possible. 
Finally we have to find a critical point of a functional defined on the real line.\\    
Solutions of $({\mathcal P}_{a,b,\l})$ are critical points
in $\Di$ of
\begin{equation*}
f_{\e}(u):=\frac {1}{2} \int_{\R^N}|x|^{-2a}|\nabla
  u|^2\,dx-\frac{\l}{2}\int_{\R^N}\frac{u^2}{|x|^{2(1+a)}}\,dx-\frac{1}{p}\int_{\R^N}\big(1+\e
  k(x)\big)\frac{u_+^p}{|x|^{bp}}\,dx,
\end{equation*}
where $u_+:=\max\{u,0\}$. For $\eps=0$ we show that $f_0$ has a one dimensional manifold of
critical points 
\[
Z_{a,b,\lambda}:= \left\{z_\mu^{a,b,\l}:=\mu^{-\frac{N-2-2a}2}z_1^{a,b,\l}\Big(\frac
x{\mu}\Big)\where \mu>0\right\},
\] 
where $z_1^{a,b,\l}$ is explicitly given in (\ref{eq:37}) below. These radial solutions
were computed for $\lambda=0$ in \cite{CatrinaWang}, the case $a=b=0$ and $-\infty
<\lambda<(N-2)^2/4$ was done by \citet{terracini}. 
The exact knowledge of the critical
manifold enables us to clarify the question of non-degeneracy.
\begin{theorem}\label{t:nondeg}
Suppose $a,b,\lambda,p$ satisfy (\ref{eq:6}). Then the critical manifold $Z_{a,b,\l}$ is non-degenerate, i.e.
\begin{align}
\label{eq:17}
T_zZ_{a,b,\l}=\ker D^2f_0(z)\quad\forall\, z\in Z_{a,b,\l},  
\end{align}
if and only if 
\begin{align}
\label{eq:18}
b\not=h_j(a,\lambda):=\frac{N}{2}\,\bigg[1+\frac{4j(N+j-1)}{(N-2-2a)^2-4\lambda}\bigg]^{-1/2}-\frac{N-2-2a}2\quad\forall\,j\in\N\setminus\{0\}.  
\end{align}
\end{theorem}
\begin{center}
 \vskip0.5truecm\noindent
 \epsfxsize=3.3in \epsfbox{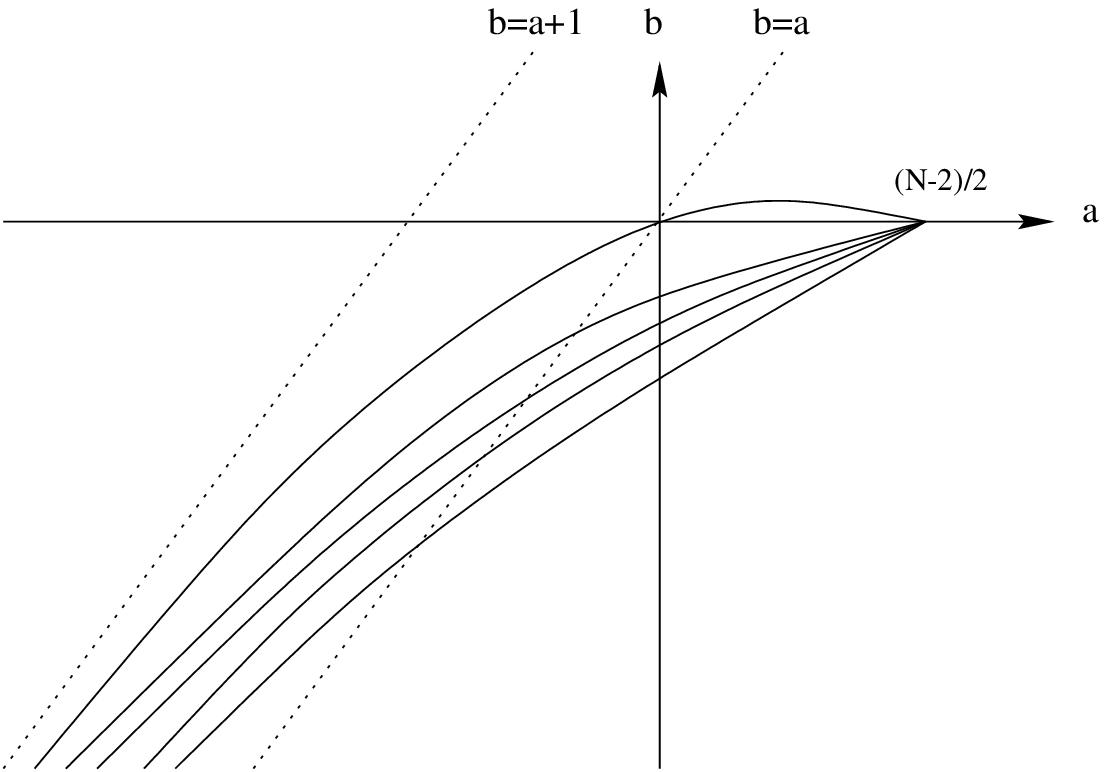}\\
{\scriptsize Figure 1 ($\lambda=0$ and $h_j(\cdot,0)$ for $j=1\dots 5$)}
\end{center}
\vskip0.5truecm\noindent
The above theorem is rather unexpected as it is explicit. It improves the
non-degeneracy results and answers an open question in \cite{AbdelPeral01}. 
Moreover, it fairly highlights the symmetry breaking
phenomenon of the unperturbed problem observed in \cite{CatrinaWang}, i.e. the existence of
non-radial minimizers of 
\begin{align}
\label{eq:28}
{{\mathcal C}_{a,b}}^{-1}:=
\inf_{u \in D_a^{1,2}(\rz^N)\backslash\{0\}}
\frac{\int|x|^{-2a}|\n u|^2}{\left(\int|x|^{-bp}|u|^p\right)^{\frac{2}{p}}}= 
\inf_{u \in D_a^{1,2}(\rz^N)  \backslash\{0\}} \frac{\|\nabla u\|^2_a}{\|u\|_{p,b}^2}.
\end{align}
In fact we improve
\cite[Thm 1.3]{CatrinaWang}, where it is shown that there are an open subset
$H\subset \rz^2$ containing $\{(a,a)\where a<0\}$, a real number $a_0\le 0$ and a function
$h:]-\infty,a_0] \to \rz$ satisfying $h(a_0)=a_0$ and $a<h(a)<a+1$ for all $a<a_0$, such that
for every $(a,b) \in H \cup \{(a,b)\in \rz^2 \where a<a_0,\,a<b<h(a)\}$ the minimizer in (\ref{eq:28})
is non-radial (see figure 2 below). We show that one may choose $a_0=0$ and $h=h_1(\cdot,0)$ and obtain, as a
consequence of Theorem \ref{t:nondeg} for $\lambda=0$,  
\begin{corollary}
\label{cor:symmbreak}
Suppose $a,b,p$ satisfy (\ref{eq:6}). If $b<h_1(a,0)$, then
${\mathcal C}_{a,b}$ in (\ref{eq:28}) is attained by a non-radially symmetric function.
\end{corollary}

\begin{center}
 \vskip0.5truecm\noindent
 \begin{tabular}{cc}
 \leavevmode 
 \epsfxsize=2.3in \epsfbox{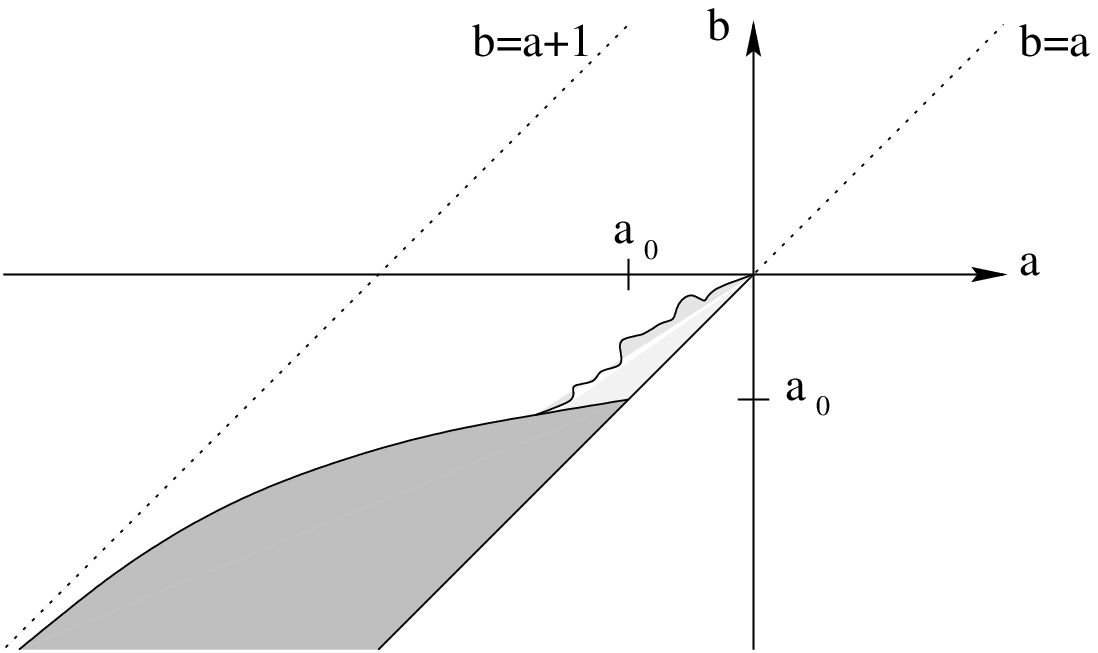} & \epsfxsize=2.3in \epsfbox{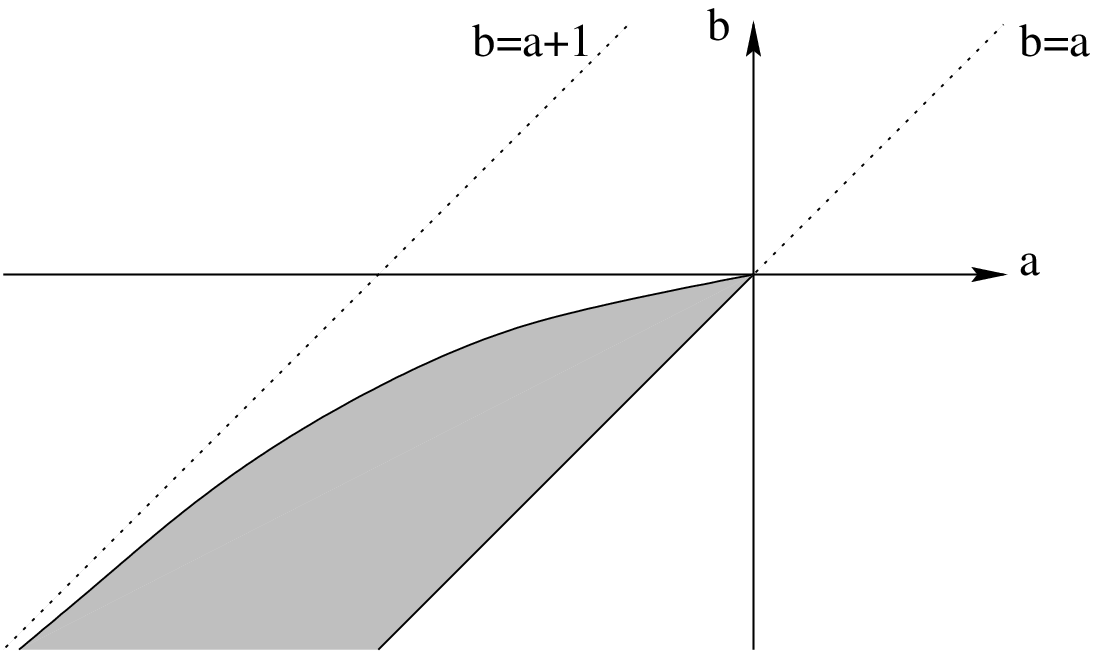}\\
 {\scriptsize region of non-radial minimizers in \cite{CatrinaWang} }&{\scriptsize region of
 non-radial minimizers given by $h_1(\cdot,0)$}\\
  \end{tabular}\\
{\scriptsize Figure 2}
\vskip0.5truecm\noindent
\end{center}

Concerning step two, the one-dimensional reduction, we follow closely the abstract scheme in
\cite{AmBa1} and construct a manifold
$Z^\eps_{a,b,\l}=\big\{z_\mu^{a,b,\l}+w(\eps,\mu)\where \mu>0\big\}$, such that any critical point of
$f_\eps$ restricted to $Z^\eps_{a,b,\l}$ is a solution to $({\mathcal P}_{a,b,\l})$. We emphasize
that in contrast to the local approach in \cite{AmBa1} we construct a manifold which is
globally diffeomorphic to the unperturbed one such that we may estimate the difference
$\|w(\eps,\mu)\|$ when $\mu \to \infty$ or $\mu \to 0$ (see also \cite{AmbGarPer01,Badiale}). 
More precisely we show
under assumption (\ref{eq:32}) below that $\|w(\eps,\mu)\|$ vanishes as $\mu \to \infty$ or
$\mu \to 0$ .\\   

We will prove the following existence results.
\begin{theorem}\label{t:main}
Suppose (\ref{eq:6}), (\ref{eq:k}), and (\ref{eq:18}) hold.
Then problem $({\mathcal P}_{a,b,\l})$ has a solution for all $|\eps|$ sufficiently small 
if
\begin{align}
\label{eq:32}
k(\infty):= \lim_{|x|\to \infty}k(x) &\text{ exists and } k(\infty)=k(0)=0.  
\end{align}
\end{theorem}
\begin{theorem}\label{t:main:2}
Assume (\ref{eq:6}),(\ref{eq:k}), (\ref{eq:18}) and
\begin{align}
\label{eq:36}
k \in C^2(\rz^N), \; |\nabla k| \in L^\infty(\rz^N) \text{ and } |D^2 k| \in L^\infty(\rz^N).  
\end{align}
Then $({\mathcal P}_{a,b,\l})$ is solvable for all small $|\eps|$ under
each of the following conditions
\begin{align}
\label{eq:33}
\limsup_{|x|\to \infty}k(x) \le k(0) &\text{ and } \laplace k(0)>0,\\
\label{eq:34}
\liminf_{|x|\to \infty}k(x) \ge k(0) &\text{ and } \laplace k(0)<0.   
\end{align}
\end{theorem} 

\begin{remark}
Our analysis of the unperturbed problem allows to consider more general perturbation, for
instance it is possible to treat equations like
\[
\begin{cases}
-\div (|x|^{-2a}\n u)-\displaystyle{\frac{\l+\eps_1 V(x)}{|x|^{2(1+a)}}}\,u=
 \big(1+\e_2 k(x)\big)\frac{u^{p-1}}{|x|^{bp}}\\
u\in D_a^{1,2}(\rz^N)   ,\quad u>0 \text{ in } \rz^N\backslash\{0\}.
\end{cases}
\]
Existence results in this direction are given by \citet{AbdelPeral01}, where the
case $a=0$ and $b=0$ is studied. We
generalize some existence results 
obtained there to arbitrary $a,b$ and $\lambda$ satisfying (\ref{eq:6}) and (\ref{eq:18}).
\end{remark}

Problem (\ref{eq:eq1}), the non-perturbative version of $({\mathcal P}_{a,b,\l})$, was 
studied by \citet{smets} in the case $a=b=0$ and $0<\lambda<(N-2)^2/4$. 
A variational minimax method combined with a careful analysis and construction of Palais-Smale
sequences shows that in dimension $N=4$ equation (\ref{eq:eq1}) has a positive solution $u\in
D_a^{1,2}(\rz^N)$ if $K\in C^2$ is positive and satisfies an analogous condition to
(\ref{eq:32}), namely $K(0)=\lim_{|x|\to\infty}K(x)$.
In our perturbative approach we need
not to impose any condition on the space dimension $N$. Theorem \ref{t:main} gives the
perspective to relax the restriction $N=4$ on the space dimension also in the nonperturbative case.

\begin{center}

{\bf Acknowledgements}

\end{center}

The authors would like to thank Prof. A. Ambrosetti for his interest
in their work and for helpful suggestions.

\section*{Preliminaries}
\citet{CatrinaWang} proved that for $b=a+1$
\[
{\mathcal C}_{a,a+1}^{-1}={\mathcal S}_{a,a+1}=\inf_{D_a^{1,2}(\rz^N) \setminus
  \{0\}}\frac{\int_{\R^N}|x|^{-2a}|\n
  u|^2}{\left(\int_{\R^N}|x|^{-2(1+a)}|u|^2\right)}=\bigg(\frac{N-2-2a}2\bigg)^2.
\]
Hence we obtain for $-\infty< \l<\Big(\frac{N-2-2a}2\Big)^2$ a norm,
equivalent to $\|\cdot \|_*$, given by
\begin{align}\label{eq:norma}
\|u\|=\left[\int_{\R^N}|x|^{-2a}|\nabla
  u|^2\,dx-\l\int_{\R^N}\frac{u^2}{|x|^{2(1+a)}}\,dx\right]^{1/2}.
\end{align}
We denote by $\Di$ the Hilbert space equipped with the scalar product induced by $\|\cdot \|$ 
\[
(u,v)=\int_{\R^N}|x|^{-2a}\nabla
  u\cdot\nabla v\,dx-\l\int_{\R^N}\frac{u\, v}{|x|^{2(1+a)}}\,dx.
\]
We will mainly work in this space. Moreover, we define by $\mathcal{C}$ the cylinder $\rz \times S^{N-1}$.
It is is shown in \cite[Prop. 2.2]{CatrinaWang} that the transformation
\begin{align}
\label{eq:20}
u(x) = |x|^{-\frac{N-2-2a}{2}}v\Big(-\ln|x|,\frac{x}{|x|}\Big)  
\end{align}
induces a Hilbert space isomorphism from $\Di$ to $H_{\lambda}^{1,2}(\mathcal{C})$, where
the scalar product in $H_{\lambda}^{1,2}(\mathcal{C})$ is defined by 
\[ (v_1,v_2)_{H_{\lambda}^{1,2}(\mathcal{C})}:= \int_{\mathcal{C}}\nabla v_1 \cdot\nabla v_2 +
\bigg(\Big(\frac{N-2-2a}{2}\Big)^2-\lambda\bigg) v_1 v_2.\]
Using the canonical identification of the Hilbert space $\Di$ with its
dual induced by the scalar-product and denoted by $\mathcal{K}$, i.e.
\[\mathcal{K}: \big(\Di\big)' \to \Di,\,  (\mathcal{K}(\varphi),u)= \varphi(u)
\quad \forall (\varphi,u) \in \big(\Di\big)'\times \Di,\]
we shall consider $f_\eps'(u)$ as an element of $\Di$ and $f_\eps''(u)$ as one of
$\mathcal{L}(\Di)$.\\
If we test $f_{\e}'(u)$ with $u_-=\max\{-u,0\}$ we get 
\begin{align*}
\big(f_{\e}'(u),u_-\big) = \int_{\R^{N}} |x|^{-2a}\nabla u \cdot\nabla u_- -\l \int_{\R^{N}} 
\frac{u u_-}{|x|^{2(1+a)}} - \int_{\R^{N}} \big(1+\e k(x)\big)\frac{u_+^{p-1}u_-}{|x|^{bp}} 
= -\|u_-\|^2  
\end{align*}
and see that any critical point of $f_{\e}$ is nonnegative. The maximum principle applied in
$\rz^N\backslash\{0\}$ shows that any nontrivial critical point is positive in that region. We
cannot expect more since the radial solutions to the unperturbed problem $(\eps=0)$ vanish at
the origin if $\lambda<0$ (see (\ref{eq:37}) below). Moreover from
standard elliptic regularity theory, solutions to $({\mathcal
  P}_{a,b,\l})$ are $C^{1,\a}(\R^N\setminus\{0\})$, $\a>0$.\\ 
The unperturbed functional $f_0$ is given by
\[
f_0(u):=\frac {1}{2} \int_{\R^N}|x|^{-2a}|\nabla
  u|^2\,dx-\frac{\l}{2}\int_{\R^N}\frac{u^2}{|x|^{2(1+a)}}\,dx-\frac{1}{p}\int_{\R^N}\frac{u_+^p}{|x|^{bp}}\,dx,\quad
  u\in \Di
\]
and we may write $f_\eps(u) = f_0(u) - \eps G(u)$, where
\begin{align}
\label{eq:29}
G(u) :=  \frac{1}{p}\int_{\R^{N}} k(x) \frac{u_+^{p}}{|x|^{bp}}. 
\end{align}

\section{The unperturbed problem}
Critical points of the unperturbed functional $f_0$ solve the equation
\begin{align}\label{eq:unperturbed}
\begin{cases}
-\div (|x|^{-2a}\n u)-\displaystyle{\frac{\l}{|x|^{2(1+a)}}}\,u=\frac
 1{|x|^{bp}}\,u^{p-1}\\
u\in \Di,\quad u>0 \text{ in } \rz^N\backslash\{0\}.
\end{cases}
\end{align}
To find all radially symmetric solutions $u$ of (\ref{eq:unperturbed}), i.e.
$u(x)=u(r)$, where $r=|x|$, we follow \cite{CatrinaWang} and note that if $u$ is radial, then equation
(\ref{eq:unperturbed}) can be written as 
\begin{align}\label{eq:radial}
-\frac{u''}{r^{2a}}-\frac{N-2a-1}{r^{2a+1}}\,u'-\frac{\l}{r^{2(a+1)}}\,u=\frac{1}{r^{bp}}\,u^{p-1}.
\end{align}
Making now the change of variable 
\begin{align}\label{eq:change1}
u(r)=r^{-\frac{N-2-2a}2}\varphi(\ln r),
\end{align}
we come to the equation
\begin{equation}\label{eq:psi}
-\varphi''+\left[\left(\frac{N-2-2a}2\right)^2-\l\right]\varphi-\varphi^{p-1}=0.
\end{equation}
All positive solutions of
(\ref{eq:psi}) in $H^{1,2}(\rz)$ are the translates of 
\begin{align*}
&\varphi_1(t)=\left[\frac{N(N-2-2a)\sqrt{(N-2-2a)^2-4\lambda}}{4(N-2(1+a-b))}\right]^{\frac{N-2(1+a-b)}{4(1+a-b)}}\cdot\\
&\quad\qquad\qquad\cdot\left(\cosh\frac{(1+a-b)\sqrt{(N-2-2a)^2-4\lambda}}{N-2(1+a-b)}\,t\right)^{-\frac{N-2(1+a-b)}{2(1+a-b)}},
\end{align*}
namely $\varphi_{\mu}(t)=\varphi_1(t-\ln \mu)$ for some $\mu>0$ (see \cite{CatrinaWang}). Consequently all radial
solutions of (\ref{eq:unperturbed}) are dilations of  
\begin{align}
\label{eq:37}
&z_1^{a,b,\l}(x)=\left[\frac{N(N-2-2a)\sqrt{(N-2-2a)^2-4\lambda}}{N-2(1+a-b)}\right]^{\frac{N-2(1+a-b)}{4(1+a-b)}}\cdot\notag\\
&\qquad\cdot\left[|x|^{\left(1-\frac
   {\sqrt{(N-2-2a)^2-4\lambda}}{N-2-2a}\right)\frac{(N-2-2a)(1+a-b)}{N-2(1+a-b)}}\Big[1+|x|^{\frac{2(1+a-b)\sqrt{(N-2-2a)^2-4\lambda}}{N-2(1+a-b)}}\Big]\right]^{-\frac{N-2(1+a-b)}{2(1+a-b)}}
\end{align}
and given by
\[
z_{\mu}^{a,b,\l}(x)=\mu^{-\frac{N-2-2a}2}z_1^{a,b,\l}\Big(\frac
x{\mu}\Big),\quad \mu>0.
\]
Using the change of coordinates in (\ref{eq:change1}), respectively (\ref{eq:20}), and the
exponential decay of $z_{\mu}^{a,b,\l}$ in these coordinates it is easy to see that the map $\mu
\mapsto z_{\mu}^{a,b,\l}$ is at least twice continuously differentiable from $(0,\infty)$ to $\Di$ and we
obtain
\begin{lemma}
\label{sec:unper:lem1}
Suppose $a,b,\lambda,p$ satisfy (\ref{eq:6}). Then the unperturbed functional $f_0$ has a one
dimensional $C^2$-manifold of critical points $Z_{a,b,\l}$ given by
$\big\{z_{\mu}^{a,b,\l}\where \mu>0\big\}$. Moreover, $Z_{a,b,\l}$ is exactly the set
of all radially symmetric, positive
solutions of (\ref{eq:unperturbed}) in $\Di$.    
\end{lemma}
In order to apply the abstract perturbation method we need to show that the manifold
$Z_{a,b,\l}$ satisfy a non-degeneracy condition. This is the content of Theorem \ref{t:nondeg}.

\begin{proof}[\bf Proof of Theorem \ref{t:nondeg}.]
The inclusion $T_{z_{\mu}^{a,b,\l}}Z_{a,b,\l}\subseteq \ker D^2
f_0(z_{\mu}^{a,b,\l})$ always holds and is a consequence of the fact that $Z_{a,b,\l}$ is a
manifold of critical points of $f_0$. Consequently, we have only to show that
$\ker D^2f_0(z_{\mu}^{a,b,\l})$ is one dimensional. Fix $u \in\ker
D^2f_0(z_{\mu}^{a,b,\l})$. The function $u$ is a solution of the linearized
problem 
\begin{align}\label{eq:linearized}
-\div (|x|^{-2a}\nabla u)
 -\frac{\l}{|x|^{2(a+1)}}u=\frac{p-1}{|x|^{bp}}(z_{\mu}^{a,b,\l})^{p-2}u.
\end{align}
We expand $u$ in spherical harmonics
\[
u(r\vartheta)=\sum_{i=0}^{\infty}\vec{v}_i(r)\vec{Y}_i(\vartheta),\quad 
r\in\R^+,\quad\vartheta\in {\mathbb S}^{N-1},
\]
where 
$\vec{v}_i(r)=\int_{{\mathbb S}^{N-1}}u(r\vartheta)\vec{Y}_i(\vartheta)\,d\vartheta$ and $\vec{Y}_i$
denotes the orthogonal $i$-th spherical harmonic jet satisfying for all $i \in \nz_0$  
\begin{align}\label{eq:harmonic}
-\Delta_{{\mathbb
    S}^{N-1}}\vec{Y}_i=i(N+i-2)\vec{Y}_i.
\end{align}
Since $u$ solves (\ref{eq:linearized}) the functions $\vec{v}_i$ satisfy for all $i\ge 0$
\[
-\frac{\vec{v}_i{''}}{r^{2a}}\,\vec{Y}_i-\frac{N-1-2a}{r^{2a+1}}\,
 \vec{v}_i{'}\,\vec{Y}_i-\frac{\vec{v}_i}{r^{2(a+1)}} \,\Delta_{\vartheta}\vec{Y}_i
 -\frac{\l}{r^{2(a+1)}}\,\vec{v}_i\,\vec{Y}_i=\frac{p-1}{r^{bp}}(z_{\mu}^{a,b,\l})^{p-2}\vec{v}_i\, 
 \vec{Y}_i 
\]
and hence, in view of (\ref{eq:harmonic}), 
\begin{align}\label{eq:vk}
-\frac{\vec{v}_i{''}}{r^{2a}}-\frac{N-1-2a}{r^{2a+1}}\, \vec{v}_i{'}+\frac{i(N+i-2)}{r^{2(a+1)}}
 \,\vec{v}_i-\frac{\l}{r^{2(a+1)}} \,\vec{v}_i=\frac{p-1}{r^{bp}}(z_{\mu}^{a,b,\l})^{p-2}\vec{v}_i.
\end{align}
Making in (\ref{eq:vk}) the transformation (\ref{eq:change1}) we obtain the equations 
\[
-\vec{\varphi}_i{''}-\beta\cosh^{-2}\big(\gamma
 (t-\ln
 \mu)\big)\vec{\varphi}_i=\left(\l-\bigg(\frac{N-2-2a}{2}\bigg)^2-i(N+i-2)\right)\vec{\varphi}_i,
\quad i \in \nz_0,
\]
where 
\[
\beta=\frac{N(N+2(1+a-b))((N-2-2a)^2-4\lambda)}{4(N-2(1+a-b))^2}\text{ and }
  \gamma=\frac{(1+a-b)\sqrt{(N-2-2a)^2-4\lambda}}{N-2(1+a-b)},
\]
which is equivalent, through the change of variable
$\zeta(s)=\varphi(s+\ln \mu)$, to
\begin{align}
\label{eq:19}
-\vec{\zeta}_i{''}-\beta\cosh^{-2}(\gamma
 s)\vec{\zeta}_i=\left(\l-\bigg(\frac{N-2-2a}{2}\bigg)^2-i(N+i-2)\right)\vec{\zeta}_i,
\quad i \in \nz_0.  
\end{align}
It is known (see \cite{olver},\cite[p. 74]{landau}) that the negative part of the spectrum of
the problem  
\[
-\zeta{''}-\beta\cosh^{-2}(\gamma s)\zeta=\nu\zeta
\]
is discrete, consists of simple eigenvalues and is given by
\[
\nu_j=-\frac{\gamma^2}{4}\left(-(1+2j)+\sqrt{1+4\beta\,\gamma^{-2}}\right)^2,\quad
j \in \nz_0,\quad 0\leq j< \frac12 \left(- 1  +\sqrt{1+4\beta\,\gamma^{-2}}\right).
\]
Thus we have for all $i \ge 0$ that zero is the only solution to (\ref{eq:19}) if and only if
\begin{align}\label{eq:seq}
A_i(a,\l)\not=B_j(a,b,\l)\text{ for all } 0\leq j<\frac{N}{2(1+a-b)},
\end{align}
where
\[
A_i(a,\l)=\l-\left(\frac{N-2-2a}{2}\right)^2-i(N+i-2)
\]
and
\[
B_j(a,b,\l)=-\frac{((N-2-2a)^2-4\lambda)(1+a-b)^2}{4(N-2(1+a-b))^2}\left[-2j+\frac
  {N}{1+a-b}\right]^2.
\]
Note that 
$A_0(a,\l)= B_1(a,b,\l)$, $A_i(a,\l)\geq A_{i+1}(a,\l)$  and  $B_j(a,b,\l)\leq
B_{j+1}(a,b,\l)$, which is shown in figure 3 below.

\vskip1truecm\noindent
\begin{center}
\leavevmode \epsfxsize=3in \epsfbox{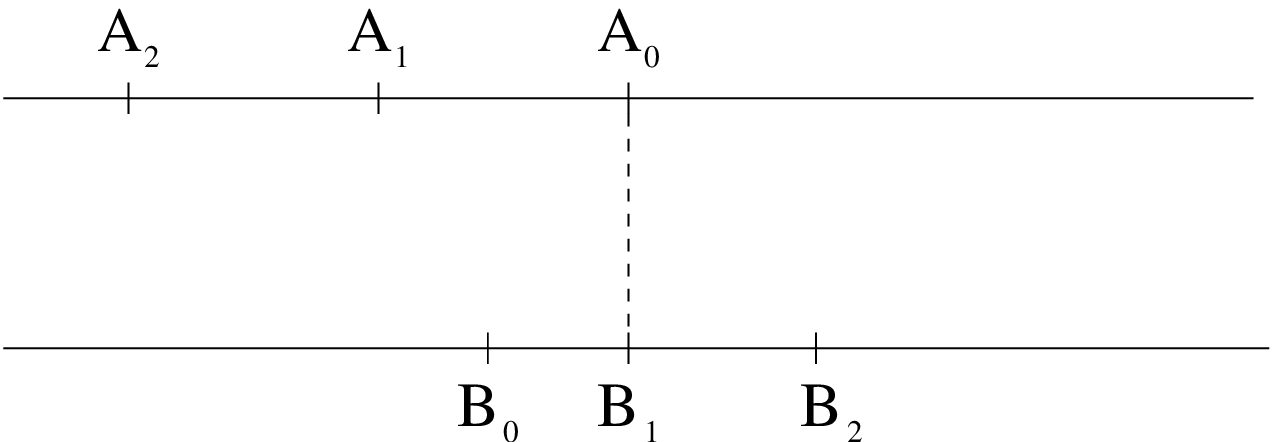}\\
{\scriptsize Figure 3}
\end{center}
\vskip1truecm\noindent

Hence (\ref{eq:seq}) is satisfied for $i\ge 1$ if and only if
$B_0(a,b,\l)\not=A_i(a,b,\l)$, which is equivalent to
$b\not=h_i(a,\l)$. On the other hand for $i=0$ equation (\ref{eq:19}) has a one
dimensional space of nonzero solutions. Hence, $\ker D^2f_0(z_{\mu}^{a,b,\l})$ is one dimensional if and only if
$b\not=h_i(a,\l)$ for any $i\geq1$, which proves the claim.  
\end{proof}

\begin{proof}[{\bf Proof of Corollary \ref{cor:symmbreak}}]
We define $I$ on $D_a^{1,2}(\rz^N)\backslash\{0\}$ by 
the right hand side of (\ref{eq:28}), i.e.
\[
I(u) := \frac{\|\nabla u\|^2_a}{\|u\|_{p,b}^2}. 
\]
$I$
is twice continuously differentiable and
\begin{align*}
(I'(u),\varphi) &= \frac{2}{\|u\|_{p,b}^2} \Big(\int_{\R^{N}}|x|^{-2a} \nabla u \nabla \varphi -
\frac{\|\nabla u\|^2_a}{\|u\|_{p,b}^p} \int_{\R^{N}} |x|^{-bp}|u|^{p-2}u\varphi \Big).
\end{align*}
Moreover, for positive critical points $u$ of $I$ a short computation leads to 
\begin{align*}
(I''(u) \varphi_1, \varphi_2) &= \frac{2}{\|u\|_{p,b}^2} 
\Big(\int_{\R^{N}}|x|^{-2a} \nabla \varphi_1 \nabla \varphi_2 
-\frac{\|\nabla u\|^2_a}{\|u\|_{p,b}^p} (p-1) \int_{\R^{N}} |x|^{-bp}u^{p-2}\varphi_1\varphi_2 \Big)\\
&\quad +(p-2) \frac{2\|\nabla u\|^2_a}{\|u\|_{p,b}^{2p+2}}
\Big(\int_{\R^{N}} |x|^{-bp}u^{p-1}\varphi_1\Big)\Big(\int_{\R^{N}} |x|^{-bp}u^{p-1}\varphi_2\Big).
\end{align*}
Obviously $I$ is constant on $Z_{a,b,0}$ and we obtain for $z_1:= z^{a,b,0}_1$ and all
$\varphi_1,\varphi_2 \in \Di$
\begin{align}
(I'(z_1),\varphi_1) &= \frac{2}{\|z_1\|_{p,b}^2} (f_0'(z_1),\varphi_1) =0, \notag\\
\label{eq:38}
(I''(z_1) \varphi_1, \varphi_2) &= 
\frac{2}{\|u\|_{p,b}^2} (f_0''(z_1)\varphi_1, \varphi_2)\notag \\
&\quad+
(p-2) \frac{2}{\|z_1\|_{p,b}^{p+2}}
\Big(\int_{\R^{N}} |x|^{-bp}z_1^{p-1}\varphi_1\Big)\Big(\int_{\R^{N}} |x|^{-bp}z_1^{p-1}\varphi_2\Big).  
\end{align}
From the proof of Theorem \ref{t:nondeg} we know that for $b<h_1(a,0)$ there exist functions
$\hat{\phi} \in D_a^{1,2}(\rz^N) $ of the form $\hat{\phi}(x)=\bar{\phi}(|x|)Y_{1}(x/ |x|)$, 
where $Y_{1}$ denotes one of the first spherical harmonics, such that
$(f_0''(z_1)\hat{\phi},\hat{\phi})<0$. By (\ref{eq:38}) we get $(I''(z_1)\hat{\phi},\hat{\phi})<0$
because the integral $\int|x|^{-bp}z_1^{p-1}\hat\phi=0$. Consequently 
${{\mathcal C}_{a,b}}^{-1}$ is strictly  
smaller than $I(z_1)=I(z^{a,b,0}_\mu)$. Since all positive radial solutions of
(\ref{eq:unperturbed}) are given by $z^{a,b,0}_\mu$ (see Lemma \ref{sec:unper:lem1}) and the
infimum in (\ref{eq:28}) is attained 
(see \cite[Thm 1.2]{CatrinaWang}) the minimizer must be non-radial.
\end{proof}
As a particular case of Theorem \ref{t:nondeg} we can state
\begin{corollary}\label{c:cor1}  
$ $
\begin{enumerate}
\item[(i)]
If\ $0<a<\frac{N-2}2$ and $0\le \lambda <\Big(\frac{N-2-2a}2\Big)^2$
then $Z_{a,b,\l}$ is non-de\-ge\-ne\-rate for any $b$ between $a$ and $a+1$.
\item[(ii)]
If $a=0$ and $0\le \lambda <\Big(\frac{N-2-2a}2\Big)^2$, then $Z_{0,b,\l}$ is degenerate if and
only if $b=\l=0$.
\end{enumerate}
\end{corollary}
\begin{remark} \label{r:Z000}
If $a=b=\l=0$, equation (\ref{eq:unperturbed}) is invariant not only
by dilations but also by translations. The manifold of
critical points is in this case $N+1$-dimensional and given by the translations and dilations of 
$z_1^{0,0,0}$. Hence the one dimensional manifold $Z_{0,0,0}$ is degenerate. 
However, the full $N+1$-dimensional critical manifold is non-degenerate in the case $a=b=\l=0$ 
(see \cite{AmbAzPer99}).
\end{remark}

\section{The finite dimensional reduction}
We follow the perturbative method
developed in \cite{AmBa1} and show that a finite dimensional reduction of our
problem is possible whenever the critical manifold is non-degenerated. 
For simplicity of notation we write $z_{\mu}$ instead
of $z_{\mu}^{a,b,\l}$ and $Z$ instead of $Z_{a,b,\l}$ if there
is no possibility of confusion.

\begin{lemma}
\label{sec:finite-dim:l1}
Suppose $a,b,\lambda,p$ satisfy (\ref{eq:6}) and $v$ is a measurable function such that 
the integral $\int |v|^{\frac{p}{p-2}}|x|^{-bp}$ is finite. Then the
operator $J_v:
D_{a,\lambda}^{1,2}(\rz^N) \to  D_{a,\lambda}^{1,2}(\rz^N)$, defined by 
\begin{align}
J_v (u) := {\mathcal K}\Big(\int_{\R^{N}}|x|^{-pb} vu \cdot\Big),  
\end{align}
is compact.
\end{lemma}
\begin{proof}
Fix a sequence $(u_n)_{n \in \nz}$
converging weakly to zero in $D_{a,\lambda}^{1,2}(\rz^N)$. 
To prove the assertion it is sufficient to show that up to a subsequence $J_v(u_n)
\to 0$ as $n \to \infty$. Using the Hilbert space isomorphism given in (\ref{eq:20}) we see
that the corresponding sequence $(v_n)_{n \in \nz}$ converges weakly to zero in
$H_\lambda^{1,2}(\mathcal{C})$. Since $(v_n)_{n \in \nz}$ converges strongly in
$L^2(\Omega)$ for all bounded domains $\Omega$ in $\mathcal{C}$,
we may extract a subsequence that converges to zero pointwise almost everywhere. Going back to
$\Di$ we may assume that this also holds for $(u_n)_{n \in \nz}$. 
By H\"older's inequality and (\ref{eq:CKN})
\begin{align*}
\|J_v(u_n)\| &\le \sup_{\|h\|_{D_{a,\lambda}^{1,2}(\rz^N)}\le 1} 
\int_{\R^N} |x|^{-pb} |v| |u_n| |h|\\
&\le  \sup_{\|h\|_{D_{a,\lambda}^{1,2}}\le 1} \Big(\int_{\R^N} |x|^{-pb} |h|^p\Big)^{1/p}
\Big(\int_{\R^N} |x|^{-pb} |v|^{\frac{p}{p-1}} |u_n|^{\frac{p}{p-1}}\Big)^{(p-1)/p}\\
&\le C \Big(\int_{\R^N} |x|^{-pb} |v|^{\frac{p}{p-1}} |u_n|^{\frac{p}{p-1}}\Big)^{(p-1)/p}.  
\end{align*}
To show that the latter integral converges to zero we use Vitali's convergence theorem given
for instance in \cite[13.38]{HewStr75}. Obviously the functions $|\cdot|^{-pb} |v|^{\frac{p}{p-1}}
|u_n|^{\frac{p}{p-1}}$ converge pointwise almost everywhere to zero. For any measurable
$\Omega\subset \rz^N$ we may estimate using H\"older's inequality
\begin{align*}
\int_{\Omega} |x|^{-pb} |v|^{\frac{p}{p-1}} |u_n|^{\frac{p}{p-1}} &\le 
\Big(\int_{\Omega} |x|^{-pb} |v|^{\frac{p}{p-2}}\Big)^{(p-2)/(p-1)} 
\Big(\int_{\Omega} |x|^{-pb} |u_n|^p \Big)^{1/(p-1)}\\
&\le C \Big(\int_{\Omega} |x|^{-pb} |v|^{\frac{p}{p-2}}\Big)^{(p-2)/(p-1)} 
\end{align*}
for some positive constant $C$. Taking $\Omega$ a set of small measure or the complement of a large ball and the use of
Vitali's convergence theorem prove the assertion.  
\end{proof}
Lemma \ref{sec:finite-dim:l1} immediately leads to
\begin{corollary}
\label{coro:compact}
For all $z \in Z$ the operator $f_0''(z): D_{a,\lambda}^{1,2}(\rz^N)
\to D_{a,\lambda}^{1,2}(\rz^N)$ may be written as $f_0''(z)= id -
J_{|z|^{p-2}}$ and is consequently a self-adjoint Fredholm operator of
index zero.    
\end{corollary}

Define for $\mu>0$ the map $U_\mu: D_{a,\lambda}^{1,2}(\rz^N) \to
D_{a,\lambda}^{1,2}(\rz^N)$ by
\begin{align*}
U_\mu(u):= \mu^{-\frac{N-2-2a}{2}} u\Big(\frac{x}{\mu}\Big).  
\end{align*}
It is easy to check that $U_\mu$ conserves the norms $\|\cdot\|$ and $\|\cdot\|_{p,b}$, 
thus for every $\mu>0$ 
\begin{align}
\label{eq:7}
(U_\mu)^{-1} &= (U_\mu)^{t} = U_{\mu^{-1}} \text{ and }
f_0 = f_0 \circ U_\mu  
\end{align}
where $(U_\mu)^{t}$ denotes the adjoint of $U_\mu$. 
Twice differentiating the identity $f_0 = f_0 \circ U_\mu$ yields for all $h_1,h_2,v \in D_{a,\lambda}^{1,2}(\rz^N)$ 
\[(f_0''(v)h_1,h_2) =
(f_0''(U_\mu(v))U_\mu(h_1),U_\mu(h_2)),\] 
that is
\begin{align}
\label{eq:1}
f_0''(v) = (U_\mu)^{-1} \circ f_0''(U_\mu(v)) \circ U_\mu \quad
\forall v \in D_{a,\lambda}^{1,2}(\rz^N). 
\end{align}
Differentiating (\ref{eq:7}) we see that $U(\mu,z):= U_\mu(z)$ maps $(0,\infty)\times Z$
into $Z$, hence 
\begin{align}
\label{eq:2}
\frac{\rand U}{\rand z}(\mu,z) = U_\mu:\: T_zZ \to T_{U_\mu(z)}Z \text{ and } 
U_\mu:\: (T_zZ)^{\perp} \to (T_{U_\mu(z)}Z)^{\perp}.  
\end{align}
If the manifold $Z$ is non-degenerated the self-adjoint Fredholm operator 
$f_0''(z_1)$ maps the space $D_{a,\lambda}^{1,2}(\rz^N)$ into
$T_{z_1}Z^\perp$ and $f_0''(z_1) \in {\mathcal{L}}(T_{z_1}Z^\perp)$ is 
invertible. Consequently, using (\ref{eq:1}) and (\ref{eq:2}), we obtain in this case 
\begin{align}
\label{eq:3}
\|(f_0''(z_1))^{-1}\|_{{\mathcal{L}}(T_{z_1}Z^\perp)} = \|(f_0''(z))^{-1}\|_{{\mathcal{L}}(T_{z}Z^\perp)} \quad \forall z \in Z.  
\end{align}

\begin{lemma}\label{l:estimates}
Suppose $a,b,p,\lambda$ satisfy (\ref{eq:6}) and (\ref{eq:k}) holds. Then there exists a constant
$C_1=C_1(\|k\|_\infty,a,b,\lambda)>0$ such that for any $\mu>0$ and
for any $w\in\Di$ 
\begin{align}
\label{eq:9}
|G(z_\mu +w)| &\le C_1\big(\||k|^{1/p} z_\mu\|_{p,b}^p + \|w\|^{p}\big)\\
\label{eq:10}
\|G'(z_\mu +w)\| &\le C_1\big(\||k|^{1/p} z_\mu\|_{p,b}^{p-1} + \|w\|^{p-1}\big) \\
\label{eq:11}
\|G''(z_\mu+w)\| &\le C_1\big(\||k|^{1/p} z_\mu\|_{p,b}^{p-2} + \|w\|^{p-2}\big) . 
\end{align}
Moreover, if $\lim_{|x| \to \infty}k(x)=: k(\infty)=0=k(0)$ then
\begin{align}
\label{eq:8}
\||k|^{1/p} z_\mu\|_{p,b} \to 0 \text{ as } \mu \to \infty \text{ or } \mu \to 0.  
\end{align}
\end{lemma}
\begin{proof}
(\ref{eq:9})-(\ref{eq:11}) are consequences of (\ref{eq:CKN}) and H\"older's inequality. We
will only show (\ref{eq:11}) as (\ref{eq:9})-(\ref{eq:10}) follow analogously.
By H\"older's inequality and (\ref{eq:CKN}) 
\begin{align*}
\|G''(z_\mu +w)\| &\le(p-1) \sup_{\|h_1\|,\|h_2\|\le 1} \int_{\R^N} \frac{|k(x)|}{|x|^{bp}}
|z_\mu+w|^{p-2}|h_1||h_2|\\ 
&\le (p-1)\||k|^{1/p}\|^2_\infty \sup_{\|h_1\|,\|h_2\|\le 1} \||k|^{1/p}
(z_\mu+w)\|_{p,b}^{p-2}\|h_1\|_{p,b}\|h_2\|_{p,b}\\
&\le c(\|k\|_\infty,a,b,\l)\,\||k|^{1/p} (z_\mu+w)\|_{p,b}^{p-2}.
\end{align*}
Using the triangle inequality and again (\ref{eq:CKN}) we obtain (\ref{eq:11}).\\
Under the additional assumption $k(0)=k(\infty)=0$ estimate
(\ref{eq:8}) follows by the dominated convergence theorem and
\begin{align*}
\int_{\R^N} \frac{|k(x)|}{|x|^{bp}} z_\mu^p 
&= \int_{\R^N} \frac{|k(\mu x)|}{|x|^{bp}} z_1^p.      
\end{align*}
\end{proof}

\begin{lemma}\label{p:implicit}
Suppose $a,b,p,\lambda$ satisfy (\ref{eq:6}) and (\ref{eq:k}) and (\ref{eq:17})
hold. 
Then there exist constants $\e_0,C>0$ and a smooth function 
\[
w=w(\mu,\e):\quad (0,+\infty) \times (-\e_0,\e_0)\ \longrightarrow\
\Di
\]
such that for any $\mu>0$ and $\e\in(-\e_0,\e_0)$
\begin{align}
\label{eq:12}
w(\mu,\e)\ \text{ is orthogonal to }\ T_{z_{\mu}}Z\\
\label{eq:13}
f_{\e}'\big(z_{\mu}+w(\mu,\e)\big)\in T_{z_{\mu}}Z\\
\label{eq:14}
\|w(\mu,\e)\|\leq C\,|\e|.
\end{align}
Moreover, if (\ref{eq:32}) holds then
\begin{align}
\label{eq:15}
\|w(\mu,\e)\|\to 0 \text{ as }\mu \to 0 \text{ or }\mu \to \infty.
\end{align}
\end{lemma}
\begin{proof}
Define $H: (0,\infty)\times D_{a,\lambda}^{1,2}(\rz^N) \times \rz \times \rz \to D_{a,\lambda}^{1,2}(\rz^N) \times \rz$
\begin{align*}
H(\mu,w,\alpha,\eps):= (f_\eps'(z_\mu+w)-\alpha \dot{\xi}_\mu,(w,\dot{\xi}_\mu)),  
\end{align*}
where $\dot{\xi}_\mu$ denotes the normalized tangent vector $\frac{d}{d\mu} z_\mu$. If
$H(\mu,w,\alpha,\eps)=(0,0)$ then $w$ satisfies (\ref{eq:12})-(\ref{eq:13}) and
$H(\mu,w,\alpha,\eps)=(0,0)$ if and only if $(w,\alpha)= F_{\mu,\eps}(w,\alpha)$, where 
\begin{align*}
F_{\mu,\eps}(w,\alpha):= -\bigg(\frac{\rand H}{\rand(w,\alpha)}(\mu,0,0,0)\bigg)^{-1}
H(\mu,w,\alpha,\eps) + (w,\alpha).
\end{align*}
We prove that $F_{\mu,\eps}(w,\alpha)$ is a contraction in
some ball $B_\rho(0)$, where we may choose the radius $\rho=\rho(\eps)>0$
independent of $z \in Z$. 
To this end we observe 
\begin{align}
\label{eq:5}
\bigg(\bigg(\frac{\rand H}{\rand(w,\alpha)}(\mu,0,0,0)\bigg)(w,\beta),
(f_0''(z_\mu)w-\beta \dot{\xi}_\mu,(w,\dot{\xi}_\mu))\bigg) =
\|f_0''(z_{\mu})w\|^2 + \beta^2+ |(w,\dot{\xi}_\mu)|^2, 
\end{align}
where
\begin{align*}
\bigg(\frac{\rand H}{\rand(w,\alpha)}(\mu,0,0,0)\bigg)(w,\beta) =
(f_0''(z_{\mu})w - \beta \dot{\xi}_\mu,(w,\dot{\xi}_\mu)).  
\end{align*}
From Corollary \ref{coro:compact} and (\ref{eq:5}) we infer that $\big(\frac{\rand H}{\rand(w,\alpha)}
(\mu,0,0,0)\big)$ is an injective Fredholm operator of index zero, hence invertible and by
(\ref{eq:3}) and (\ref{eq:5}) we obtain  
\begin{align}
\label{eq:4}
\bigg\|\bigg(\frac{\rand H}{\rand(w,\alpha)}(\mu,0,0,0)\bigg)^{-1}\bigg\| \le
\max\big(1,\|(f_0''(z_\mu))^{-1}\|\big) 
= \max\big(1,\|(f_0''(z_1))^{-1}\|\big) =: C_*.   
\end{align}
Suppose $(w,\alpha)\in B_\rho(0)$. We use (\ref{eq:1}) and (\ref{eq:4}) to see
\begin{align}
\|F_{\mu,\eps}(w,\alpha)\| &\le C_* \bigg\|\Big(H(\mu,w,\alpha,\eps) - 
\Big(\frac{\rand H}{\rand(w,\alpha)}(\mu,0,0,0)\Big)(w,\alpha)\Big)\bigg\|\notag \\
&\le C_*\|f_\eps'(z_\mu+w)-f_0''(z_\mu)w\|\notag \\
&\le C_*\int_0^1 \|f_0''(z_\mu+tw)-f_0''(z_\mu)\|\text{ dt }  \|w\| +
C_*|\eps| \|G'(z_{\mu}+w)\|\notag \\ 
&\le C_*\int_0^1 \|f_0''(z_1+tU_{\mu^{-1}}(w))-f_0''(z_1)\|\text{ dt } \|w\| +
C_*|\eps| \|G'(z_\mu+w)\|\notag\\ 
\label{eq:16}    
&\le C_*\rho \sup_{\|w\|\le \rho }\|f_0''(z_1+w)-f_0''(z_1)\| 
+ C_*|\eps| \sup_{\|w\|\le \rho }\|G'(z_\mu+w)\| .
\end{align}
Analogously we get for $(w_1,\alpha_1),(w_2,\alpha_2) \in B_\rho(0)$
\begin{align*}
\frac{\|F_{\mu,\eps}(w_1,\alpha_1)-F_{\mu,\eps}(w_2,\alpha_2)\|}{C_* \|w_1-w_2\|} &\le   
\frac{\|f_\eps'(z_\mu+w_1) - f_\eps'(z_\mu+w_2)-f_0''(z_\mu)(w_1-w_2)\|}{\|w_1-w_2\|} \\
&\le \int_0^1 \|f_\eps''(z_\mu+w_2+t(w_1-w_2))-f_0''(z_\mu)\|\di t \\
&\le \int_0^1 \|f_0''(z_\mu+w_2+t(w_1-w_2))-f_0''(z_\mu)\|
\di t \\
&\quad +|\eps| \int_0^1 \|G''(z_\mu+w_2+t(w_1-w_2))\|
\di t \\
&\le \sup_{\|w\|\le 3\rho }\|f_0''(z_1+w)-f_0''(z_1)\| + |\eps|
\sup_{\|w\|\le 3 \rho}\|G''(z_\mu+w)\|.  
\end{align*}
We may choose $\rho_0>0$ such that 
\[C_* \sup_{\|w\|\le 3\rho_0 }\|f_0''(z_1+w)-f_0''(z_1)\|<\frac{1}{2}\]
and $\eps_0>0$ such that 
\begin{align*}
2\eps_0<\Big(\sup_{z \in Z, \|w\|\le 3 \rho_0}\|G''(z+w)\|\Big)^{-1}C_*^{-1} \text{ and }
3\eps_0< \Big(\sup_{z \in Z, \|w\|\le \rho_0}\|G'(z+w)\|\Big)^{-1}C_*^{-1}\rho_0.
\end{align*}
With these choices and the above estimates it is easy to see that for every $z_\mu \in Z$ and
$|\eps|<\eps_0$ the map $F_{\mu,\eps}$ maps $B_{\rho_0}(0)$ in itself and is a contraction
there. Thus
$F_{\mu,\eps}$ has a unique fixed-point $(w(\mu,\eps),\alpha(\mu,\eps))$ in $B_{\rho_0}(0)$ and it is
a consequence of the implicit function theorem that 
$w$ and $\alpha$ are continuously differentiable.\\
From (\ref{eq:16}) we also infer that $F_{z,\eps}$ maps $B_\rho(0)$ into $B_\rho(0)$, whenever
$\rho\le\rho_0$ and
\[\rho> 2|\eps| \big(\sup_{\|w\|\le \rho}\|G'(z+w)\|\big)C_*.\]
Consequently due to the uniqueness of the fixed-point we have 
\begin{align*}
\|(w(z,\eps),\alpha(z,\eps))\| \le  3 |\eps| \big(\sup_{\|w\|\le \rho_0}\|G'(z+w)\|\big)C_*,
\end{align*}
which gives (\ref{eq:14}). Let us now prove (\ref{eq:15}). Set 
\[
\rho_{\mu}:=\min\bigg\{4\e_0C_*C_1\||k|^{1/p}z_{\mu}\|_{p,b}^{p-1},\rho_0,
    \Big(\frac1{8\e_0 C_1 C_*}\Big)^{\frac1{p-2}}\bigg\}
\]
where $C_1$ is given in Lemma \ref{l:estimates}. In view of
(\ref{eq:10}) we have that for any $|\e|<\e_0$ and $\mu>0$
\begin{align*}
2|\e|C_*\sup_{\|w\|\leq \rho_{\mu}}\|G'(z_{\mu}+w)\|&\leq
2|\e|C_*C_1\||k|^{1/p}z_{\mu}\|_{p,b}^{p-1}+2|\e|C_*C_1
\rho_{\mu}^{p-2}\rho_{\mu}.
\end{align*}
Since $\rho_{\mu}^{p-2}\le \frac1{8\e_0 C_1 C_*}$ we have, 
\begin{align*}
2|\e|C_*\sup_{\|w\|\leq \rho_{\mu}}\|G'(z_{\mu}+w)
&<2|\e|C_*C_1\||k|^{1/p}z_{\mu}\|_{p,b}^{p-1}+\frac12
\rho_{\mu}\leq\rho_{\mu},
\end{align*}
so that, by the above argument, we can conclude that $F_{\mu,\e}$ maps
$B_{\rho_{\mu}}(0)$ into $B_{\rho_{\mu}}(0)$. Consequently due to the
uniqueness of the fixed-point we have
\[
\|w(\mu,\e)\|\leq \rho_{\mu}.
\]
Since by (\ref{eq:8}) we have that $\rho_{\mu}\to 0$ for $\mu\to0$ and
for $\mu\to+\infty$, we get (\ref{eq:15}).
\end{proof}
Under the assumptions of Lemma \ref{p:implicit} we may define for $|\eps|<\eps_0$
\begin{align}
\label{eq:21}
Z_{a,b,\l}^\eps := \big\{u \in \Di \where u = z^{a,b,\l}_\mu + w(\mu,\eps),\; \mu \in (0,\infty)\big\}.  
\end{align}
Note that $Z^\eps$ is a one dimensional manifold.
\begin{lemma}
Under the assumptions of Lemma \ref{p:implicit} we may choose $\eps_0>0$ such that for every
$|\eps|<\eps_0$ the manifold $Z^\eps$ is a natural
constraint for $f_\eps$, i.e. every critical point of $f_\eps|_{Z^\eps}$ is a critical point of $f_\eps$.  
\end{lemma}
\begin{proof}
Fix $u \in Z^\eps$ such that $f_\eps|_{Z^\eps}'(u)=0$. In the
following we use a dot for the derivation with respect to $\mu$. 
Since $(\dot{z}_\mu,w(\mu,\eps))=0$ for
all $\mu >0$ we obtain
\begin{align}
\label{eq:22}
(\ddot{z}_\mu, w(\mu,\eps))+(\dot{z}_\mu,\dot{w}(\mu,\eps))=0.  
\end{align}
Moreover differentiating the identity $z_\mu= U_\sigma z_{\mu/\sigma}$
with respect to $\mu$ we obtain
\begin{align}
\label{eq:23}
\dot{z}_\sigma= \frac{1}{\sigma} U_\sigma \dot{z}_1 \text{ and } 
\ddot{z}_\sigma= \frac{1}{\sigma^2} U_\sigma \ddot{z}_1. 
\end{align}
From (\ref{eq:13}) we get that $f_\eps'(u)=c_1 \dot{z}_\mu$ for some
$\mu>0$. By (\ref{eq:22}) and (\ref{eq:23})
\begin{align*}
0 &= (f_\eps'(u),\dot{z}_\mu+\dot{w}(\mu,\eps)) = c_1
(\dot{z}_\mu,\dot{z}_\mu+\dot{w}(\mu,\eps))\\
&= c_1 \mu^{-2} \big(\|\dot{z}_1\|^2-(\ddot{z}_1, U_{\mu^{-1}}w(\mu,\eps))\big)
= c_1 \mu^{-2}\big(\|\dot{z}_1\|^2-\|\ddot{z}_1\|O(1)\eps)\big).   
\end{align*}
Finally we see that for small $\eps>0$ the number $c_1$ must be zero and
the assertion follows.  
\end{proof}
In view of the above result we end up facing a finite dimensional problem as it is enough to
find critical points of the functional $\Phi_{\e}: (0,\infty)\to \rz$ given by $f_\eps|_{Z^\eps}$.

\section{Study of $\Phi_{\e}$}
In this section we will assume that the critical manifold is non-degenerate, i.e. (\ref{eq:17}),
such that the functional $\Phi_\eps$ is defined. To find critical points of
$\Phi_{\e}=f_\eps|_{Z^\eps}$ it is convenient to introduce the 
functional $\Gamma$ given below. 
\begin{lemma}
\label{lem:gamma}
Suppose $a,b,p,\lambda$ satisfy (\ref{eq:6}) and (\ref{eq:k}) holds. Then
\begin{align}
\label{eq:24}
 \Phi_{\e}(\mu)= f_0(z_1)-\eps \Gamma(\mu)+o(\eps),  
\end{align}
where $\Gamma(\mu)=G(z_\mu)$. In particular, there is $C>0$, independent of $\mu$ and
$\eps$, such that
\begin{align}
\label{eq:35}
|\Phi_{\e}(\mu)- (f_0(z_1)-\eps \Gamma(\mu))|\le 
C \big(\|w(\e,\mu)\|^2+(1+|\e|)\|w(\e,\mu)\|^p+|\eps|\|w(\e,\mu)\| \big).    
\end{align}
Consequently, if there exist $0<\mu_1<\mu_2<\mu_3<\infty$ such that
\begin{align}
\label{eq:31}
\Gamma(\mu_2)>\max(\Gamma(\mu_1),\Gamma(\mu_3)) 
\text{ or } \Gamma(\mu_2)<\min(\Gamma(\mu_1),\Gamma(\mu_3)) 
\end{align}
then $\Phi_\eps$ will have a critical point, if $\eps>0$ is sufficiently small.
\end{lemma}
\begin{proof}
Note that for all $\mu>0$ we have $f_0(z_{\mu})=f_0(z_1)$,  
\begin{align}\label{eq:phie1}
\|z_{\mu}\|^2=\int_{\R^N}\frac{z_{\mu}^{p}}{|x|^{bp}}\text{ and }
\big(z_{\mu},w(\e,\mu)\big)=\int_{\R^N}\frac{z_{\mu}^{p-1}w(\e,\mu)}{|x|^{bp}}. 
\end{align}
From (\ref{eq:phie1}) we infer
\[
\Phi_{\e}(\mu)=\frac12\int_{\R^N}\frac{z_{\mu}^{p}}{|x|^{bp}}+
\frac12\|w(\e,\mu)\|^2+\int_{\R^N}\frac{z_{\mu}^{p-1}w(\e,\mu)}{|x|^{bp}}
-\frac 1p\int_{\R^N}\frac{(1+\e k)\big(z_{\mu}+w(\e,\mu)\big)^p}{|x|^{bp}}
\]
and
\[
f_0(z_1)=f_0(z_{\mu})=\frac 12 \|z_{\mu}\|^2-\frac
1p\int_{\R^N}\frac{z_{\mu}^{p}}{|x|^{bp}}=\bigg(\frac 12-\frac
1p\bigg) \int_{\R^N}\frac{z_{\mu}^{p}}{|x|^{bp}}.
\]
Hence
\begin{align}
\label{eq:27}
\Phi_{\e}(\mu)=f_0(z_1)-\eps \Gamma(\mu) +\frac12\|w(\e,\mu)\|^2-\frac 1p H_{\e}(\mu),  
\end{align}
where
\[
H_{\e}(\mu)=\int_{\R^N}
\frac{\big(z_{\mu}+w(\e,\mu)\big)^p-z_{\mu}^p-p\,z_{\mu}^{p-1}w(\e,\mu) 
+ \eps k \big((z_{\mu}+w(\e,\mu)\big)^p-z_{\mu}^p\big)}{|x|^{bp}}.
\]
Using the inequality
\begin{equation*}
(z+w)^{s-1}-z^{s-1}-(p-1)z^{s-2}w\leq
\begin{cases}
C (z^{s-3}w^2+w^{s-1})&\mbox{if}\ s\geq 3\\
C\ w^{s-1}&\mbox{if}\ 2<s<3,
\end{cases}
\end{equation*}
where $C=C(s)>0$, with $s=p+1$ and H\"older's inequality we have for
some $c_2,c_3>0$
\begin{align*}
|H_{\e}(\mu)|&\leq
\int_{\R^N}\frac{\big|\big(z_{\mu}+w(\e,\mu)\big)^p-z_{\mu}^p
-p\,z_{\mu}^{p-1}w(\e,\mu)\big|}{|x|^{bp}}
+|\e|\int_{\R^N}\frac{|k|\, \big((z_{\mu}+w(\e,\mu))^p -z_\mu^p\big)}{|x|^{bp}}\\
&\leq c_2\left[\int_{\R^N}\frac{z_{\mu}^{p-2}w^2(\e,\mu)}{|x|^{bp} }
+\int_{\R^N}\frac{|w(\e,\mu)|^p }{|x|^{bp}}
+|\e|\int_{\R^N}\frac{z_{\mu}^{p-1}|w(\eps,\mu)|}{|x|^{bp}}
+|\e|\int_{\R^N}\frac{|w(\e,\mu)|^p }{|x|^{bp}}\right]\\
&\leq c_3\left[\|w(\e,\mu)\|^2+(1+|\e|)\|w(\e,\mu)\|^p+|\e|\|w(\e,\mu)\|\right]
\end{align*}
and the claim follows.
\end{proof}
Although it is convenient to study only the reduced functional $\Gamma$ instead of $\Phi_\eps$,
it may lead in some cases to a loss of information, i.e. $\Gamma$ may be constant even if $k$
is a non-constant function. This is due to the fact that the critical manifold consists of
radially symmetric functions. Thus $\Gamma$ is constant for every $k$ that has constant
mean-value over spheres, i.e. 
\[ \frac{1}{r^{N-1}}\int_{\rand B_r(0)} k(x) \di S(x)\equiv \text{const } \quad \forall r>0 \text.\]  
In this case we have to study the functional $\Phi_{\e}(\mu)$ directly.
\begin{proof}[\bf Proof of Theorem \ref{t:main}]
By (\ref{eq:32}), (\ref{eq:8}), (\ref{eq:15}) and (\ref{eq:35}) 
\[
\lim_{\mu\to 0^+}\Phi_{\e}(\mu)=\lim_{\mu\to
  +\infty}\Phi_{\e}(\mu)=f_0(z_1).
\]
Hence, either the functional $\Phi_{\e}\equiv f_0(z_1)$, and we obtain infinitely many critical
points, or $\Phi_{\e}\not\equiv f_0(z_1)$ and $\Phi_{\e}$ has at least a 
global maximum or minimum. In any case $\Phi_{\e}$ has a critical point that
provides a solution of $({\mathcal P}_{a,b,\l})$.
\end{proof}
The next lemma shows that it is possible (and convenient) to extend the $C^2-$ functional
$\Gamma$ by continuity to $\mu=0$. The proof of this fact is analogous to the one in
\cite[Lem. 3.4]{AmbAzPer99} and we omit it here.  
\begin{lemma}
\label{lem:gamma2}
Under the assumptions of Lemma \ref{lem:gamma} 
\begin{align}
\label{eq:25}
&\Gamma(0):= \lim_{\mu \to 0} \Gamma(\mu) = k(0) \frac{1}{p} \|z_1\|^p_{p,b} \quad\text{ and } \\
\label{eq:26}
&\frac{1}{p}\liminf_{|x|\to \infty}k(x)\|z_1\|^p_{p,b} \le \liminf_{\mu \to \infty} \Gamma(\mu) \le 
\limsup_{\mu \to \infty} \Gamma(\mu) \le  \frac{1}{p}\limsup_{|x|\to \infty}k(x)\|z_1\|^p_{p,b}.  
\end{align}
If, moreover, (\ref{eq:36}) holds we obtain
\begin{align}
\label{eq:30}
\Gamma'(0)=0 \text{ and } \Gamma''(0)=\frac{\laplace k(0)}{Np} \int |x|^2 \frac{z_1(x)^p}{|x|^{bp}}.  
\end{align}
\end{lemma}
\begin{proof}[\bf Proof of Theorem \ref{t:main:2}]
To see that assumptions (\ref{eq:33}) and (\ref{eq:34}) give rise to a
critical point we use the functional $\Gamma$. Condition (\ref{eq:33})
and Lemma \ref{lem:gamma2} imply that $\Gamma$ 
has a global maximum strictly bigger than $\Gamma(0)$ and $\limsup_{\mu \to \infty}
\Gamma(\mu)$. Consequently $\Phi_{\e}$ has a critical point in view of Lemma
\ref{lem:gamma}. The same reasoning yields a critical point under condition (\ref{eq:34}).   
\end{proof}

\bibliographystyle{adinat}
\bibliography{schneiderfelli}

\end{document}